\documentclass[a4paper,12pt]{article}

\usepackage{amssymb}
\usepackage{graphicx}
\usepackage{dcolumn}
\usepackage{amsmath,amssymb,mathrsfs}

\begin{document}

\title{Sophomore's dream  function:  asymptotics,  complex plane behavior and relation to the error function}

\author{V. Yu. Irkhin}
\maketitle
{\it Institute of Metal Physics, 620108 Ekaterinburg, Russia}
\begin{abstract}
Sophomore's dream sum $S=\sum_{n=1}^\infty n^{-n}$ is extended to the function $f(t,a)=t\int_{0}^{1}(ax)^{-tx}dx$ with $f(1,1)=S$.  Asymptotic behavior for a large  $|t|$ is obtained, which is exponential for $t>0$ and $t<0,a>1$, and inverse-logarithmic for $t<0,a<1$. An advanced approximation includes  a half-derivative of the exponent and is expressed in terms of the  error function. This approach  provides excellent interpolation  description in the complex plane. The function $f(t,a)$ demonstrates  for $a>1$ oscillating behavior  along the imaginary axis with slowly increasing amplitude and the period of $2\pi iea$,   modulation by high-frequency oscillations being present.  Also,  $f(t,a)$  has non-trivial zeros in the left complex half-plane with Im~$t_n \simeq 2(n-1/8)\pi e/a$ for $a\geq 1$. The results  describe analytical integration of the function $x^{tx}$.
\end{abstract}

\section{Introduction}
One of most beautiful mathematical results is provided by the identities
discovered by Johann Bernoulli in 1697 \cite{552,553}:%
\begin{eqnarray}
\sum_{n=1}^{\infty }\frac{1}{n^{n}} &=&\int_{0}^{1}x^{-x}dx \simeq 1.2913, \label{s1} \\
\sum_{n=1}^{\infty }\frac{(-1)^{n-1}}{n^{n}} &=&\int_{0}^{1}x^{x}dx \simeq 0.7834 .\label{s2}
\end{eqnarray}%
They can be proven by expanding $x^{x}=\exp (x\ln x)$ in the power series
and subsequent integration. Later they were named the ``sophomore's dream''
owing to their counterintuitive simplicity, generalized and investigated in a number of works (see, e.g., \cite{Pi,Pi1}).

The sophomore's dream derivation can be generalized to the function
\begin{equation}
f(t)=\sum_{n=1}^{\infty }\frac{t^{n}}{n^{n}}  \label{f}
\end{equation}%
which satisfies
\begin{equation}
f(t)=t\int_{0}^{1}x^{-tx}dx.  \label{s3}
\end{equation}%
In particular, for $t=1$ we come to Eq.(\ref{s1}), and for $t=-1$  to Eq.(\ref{s2}).
Indeed, expanding  $x^{-tx}=\exp (-tx\ln x)$ we obtain the integral
\begin{equation}
	\int_{0}^{1}x^{-tx}dx=\int_{0}^{1} \sum_{n=0}^{\infty }\frac{(-tx)^{n}\ln^n(x) }{n!} dx.
	  \label{s31}
\end{equation}%
Performing the substitution $x = \exp (-\frac{z}{n+1})$ we have
\begin{equation}
	\int_{0}^{1}x^{-tx}dx
	=  \sum_{n=0}^{\infty }	\frac{t^{n} }{n!} (n+1)^{-(n+1)}\Gamma (n+1)
	=\sum_{n=1}^{\infty }
	\frac{t^{n-1} }{n^{n}}  \label{s311}
\end{equation}%
where
 $$\Gamma (n)=\int_{0}^{\infty}z^{n-1} \exp(-z) dz = (n-1)!
$$ is  Euler's gamma function.

A more general integral
\begin{equation}
	I(a,b,c,d)=\int_{0}^{1}(ax+b)^{cx+d}dx  \label{s13}
\end{equation}%
was treated in Ref. \cite{Pi}. It is transformed into $n$-sum in a similar way, but contains, unlike Eq.(\ref{s311}), the incomplete upper gamma function
\begin{equation}
	\Gamma (y,z)=\int_{z}^{\infty}x^{y-1} \exp(-x) dx. \label{s131}
\end{equation}%
In particular, 
\begin{equation}
I(a,0,c,d)=\frac{1}{a}\sum_{n=0}^{\infty}\frac{(-c/a)^{n}}{n!(n+d+1)^{n+1}}\Gamma(n+1,-(n+d+1)\ln a) \label{s161}
\end{equation}%
and for $a=1$
\begin{equation}
	I(1,0,-t,d)= \sum_{n=1}^{\infty }
	\frac{t^{n-1}}{(n+d)^{n}}
	\label{s16}.
	\end{equation}%
which generalizes Eqs. (\ref{f}), (\ref{s3}).

Here we treat the function
\begin{equation}
	f(t,a)=t\int_{0}^{1}(ax)^{-tx}dx =t I(a,0,-t,0)	  \label{sa}.
\end{equation}%
We have
\begin{equation}
	\Gamma(n+1,-(n+1)\ln a)=n! - \int_{0}^{(n+1)\ln(1/a)}\exp ({n\ln x-x)}dx.
	\label{gam}
\end{equation}
For large $n$ (corresponding to large $t>0$) and $a>1/e$ the integral in (\ref{gam}) is very small and $f(t,a)$ is a function of $t/a$ only:
\begin{equation}
	f(t,a)\simeq \sum_{n=1}^{\infty }
	\frac{(t/a)^{n}}{n^{n}}	  \label{sa6}.
\end{equation}%
For $a<1/e$, the situation becomes more complicated since the integral contains a saddle point $x=n$.
Thus we have a crossover point $a=1/e$, which will be discussed below.
The case $t<0$ also requires a special consideration.

It is instructive to seek for analogies of the function (\ref{sa})  with the usual exponential function
\begin{equation*}
\exp t=\sum_{n=0}^{\infty }\frac{t^{n}}{n!},
\end{equation*}%
which has wide physical applications. We perform analytical and numerical investigation of  properties of the functions $f(t)=f(t,1)$ and $f(t,a)$ both on the real axis and in the complex plane.
To this end we establish connection with  well-known special functions. By the way, the problem of analytical integration of the function $x^x$ becomes more clear.

\section{Half-derivative representation  in terms of special functions}
First we treat Eq.(\ref{f}) by  using in the sum the Stirling formula
\begin{equation}
n!\ =(2\pi n)^{1/2}\frac{n^{n}}{e^{n}},
\end{equation}%
which works well even at small $n$, to obtain%
\begin{equation}
f(t)=(2\pi )^{1/2}\sum_{n=1}^{\infty }n^{1/2}\frac{(t/e)^{n}}{n!}.
\end{equation}%
Making the approximation $\Gamma (n+1)/\Gamma (n+1/2) \simeq n^{1/2}$ we may transform
\begin{equation}
f(t)\simeq \left( \frac{2\pi t}{e}\right) ^{1/2}\sum_{n=1}^{\infty }\frac{%
\Gamma (n+1)}{\Gamma (n+1/2)}\frac{(t/e)^{n-1/2}}{n!}=\left( \frac{2\pi t}{e}%
\right) ^{1/2}F^{(1/2)}(\frac{t}{e})  \label{half}
\end{equation}%
 where
\begin{eqnarray}
F^{(1/2)}(x) =\frac{\partial ^{1/2}}{\partial x^{1/2}}\exp x-\frac{1}{%
\sqrt{\pi x}}
\end{eqnarray}%
is nothing but the half-derivative of the exponential function with the
singular part (half-derivative of unity)  being excluded.
We will see that this approximation provides a good  accuracy even at small $t$ where small $n$ work in the sum.

The half-derivative  is defined by the integral  \cite{frac}
\begin{equation}
F^{(1/2)}(x) = \pi^{-1/2}\int_{0}^{x}\frac{\exp y}{\sqrt{x-y}}dy
	=\pi^{-1/2}\gamma (\frac{1}{2},x)\exp x \label{haf}
\end{equation}%
where $\gamma (\alpha ,x)$ is the lower incomplete gamma function,
which has the series expansion
\begin{equation}
	\gamma (1/2,x)=(\pi x)^{1/2}\exp (-x) \sum_{n=1}^{\infty }\frac{x^{n-1}}{\Gamma(n+1/2)}  \label{fg}
\end{equation}%
and
\begin{equation}
	\frac{\partial ^{1/2}x^{n}}{\partial x^{1/2}}=\frac{\Gamma (n+1)}{\Gamma
		(n+1/2)}x^{n-1/2},~n\geq 0.
\end{equation}%
Note that the singular part corresponding to the $n=0$ term of the sum in (\ref{half}) may also be related to introducing  $n=0$  into (\ref{f}).

The integral (\ref{haf}) can be transformed as
\begin{equation}
	F^{(1/2)}(x) = {\rm erf}(\sqrt x)\exp x \label{haf11}
\end{equation}%
where
\begin{equation}
	{\rm erf}(x) =2 \pi^{-1/2}\int_{0}^{{x}}{\exp (-z^2)}dz  \label{erf}
\end{equation}%
is the error function.

At large $x,~F^{(1/2)}(x)=\exp x$ and we come  to the exponential asymptotics,  which will be discussed in detail in the next Section.
This result can be generalized to $a\neq 1$ (for $a> 1/e$), as follows from (\ref{sa6}).

At the same time, Eq. (\ref{half}) provides a
 correct $t$-linear behavior at small $t$, corresponding to the
$n=1$ term. Since ${\rm erf}(x) \simeq 2x/\sqrt \pi$ at small $x$, we have a minor error owing to the factor of $2^{3/2}/e = 1.04$.


\section{Asymptotics and interpolation}

\subsection{$t>0$}
Now we treat the  integral (\ref{s3}). Changing the
variables $x=e^{-y}$, it can be rewritten as%
\begin{equation}
	f(t,a)=t\int_{0}^{\infty }\exp(te^{-y}(y-\ln a)-y)dy. \label{y}
\end{equation}%
For large $t$, this exponential  integral can be estimated by the Laplace
method, i.e., by expanding the exponent near the  point where its derivative vanishes, $y_0\simeq 1+\ln a$. Putting $y=y_0+\delta $ and expanding in $\delta $
we get
\begin{equation*}
	te^{-y}(y - \ln a )\simeq (t/ea)(1-\delta ^{2}/2).
\end{equation*}%
Evaluating the integral in the infinite limits
\begin{equation}
	\int_{-\infty }^{\infty }\exp \left( - \frac{t\delta ^{2}}{2ea}\right) d\delta
	=\left( \frac{2\pi ea}{t}\right) ^{1/2} \label{int}
\end{equation}%
we obtain the exponential asymptotics required%
\begin{equation}
	f(t,a)=\left( \frac{2\pi t}{e a}\right) ^{1/2}\exp \frac{t}{e a}.
	\label{f16}
\end{equation}%

To obtain corrections of next order in $1/t$, we can  expand the exponent in
\begin{equation}
	f(t,a)=\frac{t}{e}\int_{-\infty}^{\infty }\exp\left(\frac{t}{e} e^{-\delta}(1+\delta-\ln a)-y\right)e^{-\delta} d\delta \label{y1}
\end{equation}%
up to $\delta^4$ and  the pre-exponential factor $e^{-\delta}$ up to $\delta^2$. Further on, we pass from $\delta$ to a new variable $s$ defined by the condition that the series in the exponent equals to $-s^2/2$. Collecting all the terms, after cumbersome transformations we express $\delta$ in terms of $s$ and derive to second order  for $a=1$
\begin{equation}
	f(t)=\frac{t}{e}\int_{-\infty}^{\infty }\exp\left(-\frac{t}{e}s^2\right)\left(1-s^2/24 \right) ds. \label{y2}
\end{equation}%
Performing integration we have
\begin{equation}
	f(t)=\sqrt{\frac{2\pi t}{e}}\left(1-\frac{e}{24t}\right)\exp  \frac{t}{e}.  \label{f3}
\end{equation}

Numerical calculations demonstrate that the  large-$t$ behavior is described to a very high accuracy by the equation%
\begin{equation}
	f(t)=\sqrt{\frac{2\pi (t-1/4)}{e}}\exp \frac{t}{e},  \label{f2}
\end{equation}
which is in agreement with (\ref{f2}).
This approximation yields also a good estimate for the integral (1), $f(1)=1.31666$.



\begin{figure}
	\includegraphics[width=0.49\textwidth]{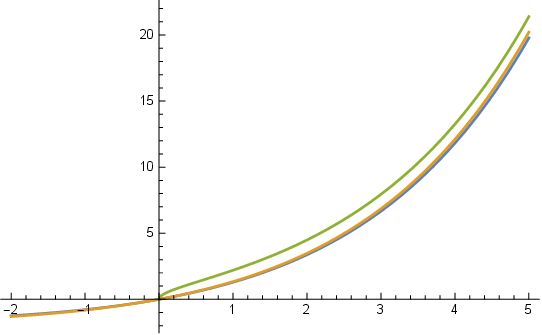}
		\includegraphics[width=0.49\textwidth]{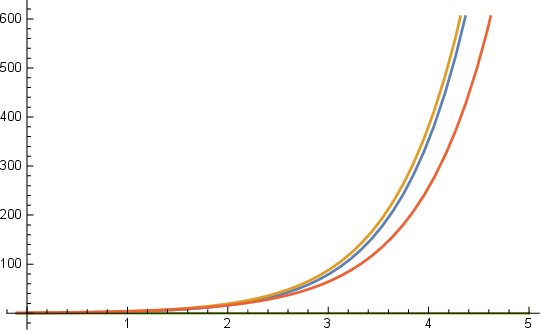}
	\caption{Left panel: Plots of  the exact $f(t)=f(t,1)$  (blue) and of the  approximations for $f(t)$ according to Eqs. (\ref{half}) (green) and (\ref{77}) (orange). The blue and orange lines are almost indistinguishable.
	Right panel:	Plots of  the exact $f(t,0.25)$  (blue) and of the  approximation (\ref{f16b}) (orange) and $f(t,a) = a^{-t}$  (red).
	}
	\label{1}       
\end{figure}

Since at large $t$ the integral (\ref{int}) depends very weakly on the integration limits, we may construct an interpolation description by choosing them from the requirement of the correct small-$t$ behavior.
Similar to (\ref{haf11}), this condition is satisfied (again with a factor of 1.04) by putting these limits to be $\mp t/(2^{3/2}e)$. Then we have
\begin{eqnarray}
	f(t,a) =\frac {\phi(ta/e)} {a}	\exp \left(\frac{t}{e a}\right), \label{77}
	\\
	\phi(x)=\left( \frac{2\pi x}{e}\right) ^{1/2} {\rm erf}(x ^{1/2})
	\label{ha11}
\end{eqnarray}%

A numerical comparison of the exact function $f(t)$ and the approximations (\ref{f16})  and (\ref{77})  is presented in  Fig.~1. One can see that the approximation (\ref{77}) works much better. It has also the advantage to work for negative $t$.

With decreasing $a$, the accuracy of the  asymptotics  (\ref{f16}) becomes lower,  and this asymptotics will not work at small $a$, since now $y_0= 1+\ln a$ is beyond the integration interval.
At the critical point $a=1/e$ the approximation (\ref{f16}) is still satisfactory, but the integration goes only at positive $\delta$, so that we have to introduce the factor of 1/2:
\begin{equation}
	f(t,1/e) = \left( {\pi t/2}\right) ^{1/2}\exp {t}.
	\label{sae}
\end{equation}%
To obtain a correct exponential dependence for $a<1/e$, we we can use  the expansion directly in the integral (\ref{sa}) near the point $x=1$, i.e., $x=1-\delta$, to obtain
\begin{equation}
	-tx(\ln x + \ln a )\simeq
		-t(\ln a +(\delta -1 -\ln a)^{2}/2 - (1 +\ln a)^2/2). \label{x2}
\end{equation}%
Replacing $\delta \rightarrow \delta'= \delta -1 -\ln a$ and performing integration in the appropriate limits  we derive
\begin{equation}
	f(t,a)=(\pi t/2)^{1/2} a^{-t} {\rm  erfc}[-(1 +\ln a)\sqrt{t/2}]  \exp (-t(1 +\ln a)^2/2)
	\label{f16b}
\end{equation}
where
\begin{equation*}
	{\rm  erfc}(x)=2 \pi^{-1/2}\int_{x}^{{\infty}}{\exp (-z^2)}dz  = 1- {\rm erf}(x)
\end{equation*}%
is the complimentary error function.

The accuracy of this result, being about 10\% for $a\simeq 0.2$, increases with decreasing $a$ since relevant region near $x=1$ becomes narrower and the expansion works better.
However, the accuracy of the simplest approximation $f(t,a) \simeq a^{-t}$ (see Fig.~1, right panel) is considerably worse even at large $t$ and especially at not small $a$. Indeed, using the asymptotic expansion
\begin{equation}
{\rm  erfc}(x)=  \frac{\exp (-x^2)}{ \pi^{1/2} x}\left( 1 + \sum_{n=1}^\infty (-1)^n\frac{(2n-1)!!}{(2x^2)^n}  \right)
	\label{fasy}
\end{equation}
we derive
\begin{equation}
	f(t,a)\simeq -(1 +\ln a)^{-1} a^{-t}[1-1/(2(1 +\ln a)^2 t)+...].
	\label{faaas}
\end{equation}
which also works satisfactorily only at small $a$.

For completeness, we consider briefly the behavior of the function
\begin{equation}
	\tilde{f}(t,a)=t\int_{1}^{\infty}(ax)^{-tx}dx	  \label{sad}.
\end{equation}%
We obtain similar to (\ref{f16b}) in the saddle-point approach
\begin{equation}
	\tilde{f}(t,a)=(\pi t/2)^{1/2} a^{-t} {\rm  erfc}[(1 +\ln a)\sqrt{t/2}]  \exp (-t(1 +\ln a)^2/2) 	  \label{sadd}.
\end{equation}%
This approximation works very well even for moderate $t$ at $a>1/e$, but at $a<1/e$ the accuracy becomes lower since the saddle point is beyond the integration region.
Thus $\tilde{f}(t,a)$ exponentially decreases for $a>1$, exponentially increases for $a<1$ with $t$, and $\tilde{f}(t\rightarrow \infty,1)\rightarrow 1$. Again,   the simplest approximation $\tilde{f}(t,a) \simeq a^{-t}$ is  very crude.

\subsection{$t<0$}
For negative $t$, the approximation (\ref{77}) works only at rather moderate $|t|$. Indeed, the point $y_0$ is now a maximum rather than a minimum of the exponent. However, for $a>1$  we can again use the expansion (\ref{x2}) to obtain
\begin{equation}
	f_1(t,a)=(-\pi t/2)^{1/2} a^{-t}  {\rm  erfi}[-(1 +\ln a)\sqrt{-t/2} ] \exp (t(1 +\ln a)^2/2)
	\label{faaa}
\end{equation}
where
\begin{equation*}
{\rm  erfi}(x)=2 \pi^{-1/2}\int_{0}^{{x}}{\exp (z^2)}dz = -i{\rm  erf}(ix)
\end{equation*}
is the imaginary error function (which is in fact real).
Thus $|f(t,a)|$ exponentially increases with $|t|$.
The result (\ref{faaa}) works very well  for $a>1, t<0$
since only a narrow region near $x=1$ works in the integral.


For $a<1$, the situation changes: now the contribution (\ref{faaa}) is exponentially small for large $|t|$. However, the true behavior is different: in fact, the integral (\ref{y}) is determined by the vicinity of the saddle point which is defined by the equation
\begin{equation}
	\frac{d}{dy}(te^{-y}(y-\ln a)-y)=te^{-y}(1+\ln a -y)-1=0.
		\label{eqs}
\end{equation}
Solving   this equation by iterations, we can approximate the integral by putting $y \simeq  \ln(-t)+\ln (\ln (-t/a)-1)...$ in the brackets of the exponent,
\begin{equation}
	f_2(t,a)\simeq t\int_{0}^{\infty} \exp ( te^{-y} [\ln(-t/a)+\ln(\ln (-t/a)-1) ]-y )dy. \label{ys}
\end{equation}%
Changing back the variables $x=e^{-y}$ we can perform integration to obtain
\begin{equation}
	f_2(t,a)= -1/[\ln(-t/a)+\ln (\ln (-t/a)-1)].
	\label{fas}
\end{equation}
This approximation works very well at large $|t|$ for any $a<1$ (note that the $\ln \ln$-correction enables one to obtain  much better accuracy than the simple logarithmic approximation). However, the region of validity of this asymptotics increases with decreasing $a$.
To make better interpolation, we can take into account the contribution (\ref{faaa}) which corrects the behavior at moderate $|t|$, i.e., put
\begin{equation}
	f(t,a)= f_1(t,a)+f_2(t,a).
	\label{fat}
\end{equation}
Vice versa, the contribution (\ref{fas}) improves somewhat the behavior at moderate $|t|$ for $a>1$ (Fig.~2).

\begin{figure}
	\includegraphics[width=0.49\textwidth]{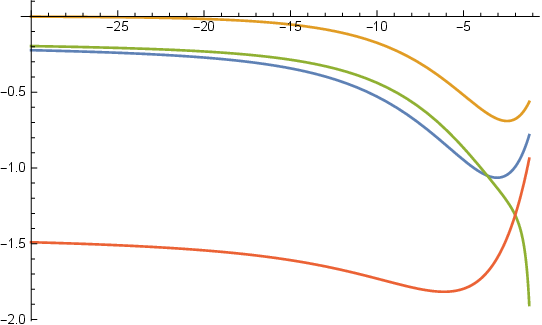}
	\includegraphics[width=0.49\textwidth]{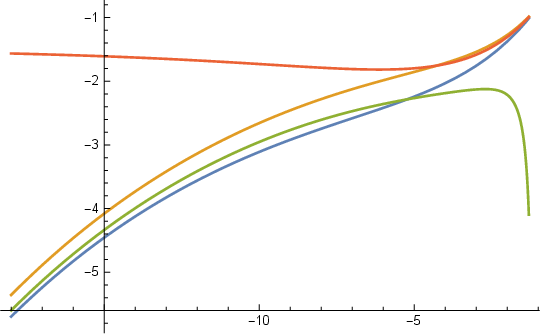}
	\caption{Plots of  the  exact $f(t,a)$  (blue) and of the  approximations according to Eqs.~(\ref{fat}) (green),  (\ref{faaa}) (orange), and (\ref{77}) (red) for $a=0.8$ (left panel) and $a=1.1$ (right panel). Note that the singular behavior of   $f_2(t,a)$ at small $|t|$ can be eliminated by exact solving equation (\ref{eqs}).
	}
	\label{22}       
\end{figure}

\begin{figure}
	\includegraphics[width=0.69\textwidth]{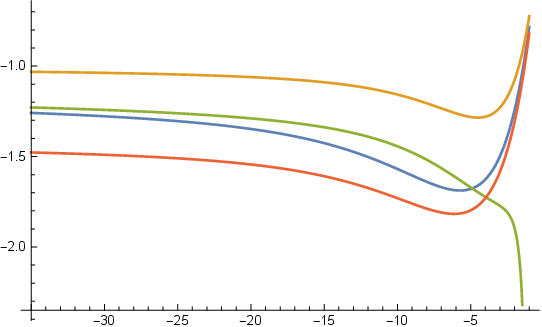}
	\caption{Plots of  the  exact $f(t)$  (blue) and of the  approximations according to Eqs.~(\ref{fat}) (green),  (\ref{faaa}) (orange), and (\ref{77}) (red).
	}
	\label{33}       
\end{figure}

Thus $f(t,a < 1)$ for large $-t>0$  is a weakly $t$-dependent function which takes  nearly constant value about  $- 0.1$ in a wide interval of  $-t$ about 10-1000. At the same time, for smaller $-t$, the approximation (\ref{77})  can be used.
We will see that such a picture  is important for the  behavior of $f(t,a)$ in the  complex plane.

For $a=1$ we have a marginal situation: $f(t\rightarrow  -\infty,1) \rightarrow -1$. This result is verified by the approximation (\ref{faaa}). At the sane time, the large-$|t|$ dependence is different from the behavior of the exact $f(t,1)$: approaching to -1 is too rapid.
Again, we can improve the interpolation  by taking into account both the contributions (Fig.~3).
On the other hand, the approximation (\ref{77}) is not  satisfactory for large $|t|$  yielding $f(t\rightarrow  -\infty,1) \rightarrow -\sqrt 2$.

The slow behavior  $f(t\rightarrow  -\infty) \rightarrow -1$ is an important property of the function $f(t)$ defined by (\ref{f}). This demonstrates again some similarity to the usual exponential function.

\section{The behavior in the complex plane}

The asymptotics (\ref{f16}) can be applied in the complex plane ($t \rightarrow z = x+iy$). One can see that oscillations with a slowly increasing amplitude occur on the imaginary axis,  but the agreement with the numerical results for the real part is not perfect even at large $|z|$, since it has, besides oscillating part, a constant contribution (Fig.~4). This shift is connected with the above-discussed behavior at real negative $t$, which is satisfactorily described by the approximation (\ref{77}).

\begin{figure}
	\includegraphics[width=0.49\textwidth]{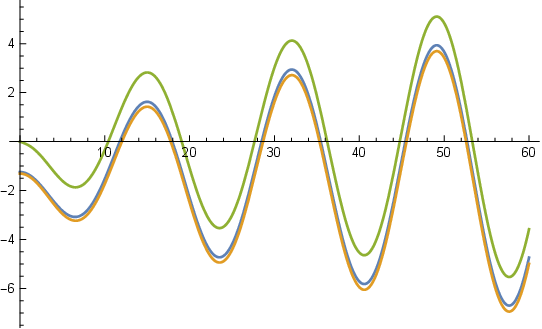}
		\includegraphics[width=0.49\textwidth]{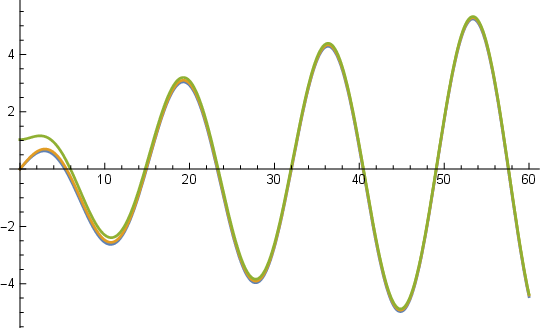}
	\caption{Plots of  the real (left panel)  and imaginary (right panel) for the exact $f(iy-2)$  (blue) and for the  approximations according to Eqs.~(\ref{half}) (green) and (\ref{77}) (orange).
	}
	\label{2}       
\end{figure}

\begin{figure}
	\includegraphics[width=0.49\textwidth]{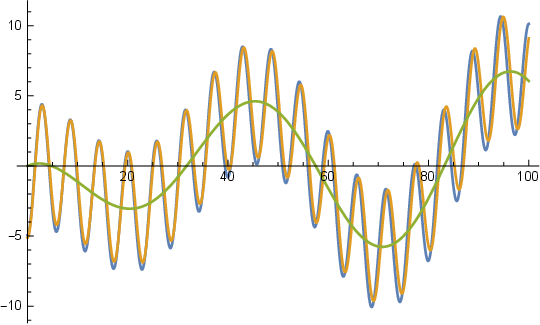}
		\includegraphics[width=0.49\textwidth]{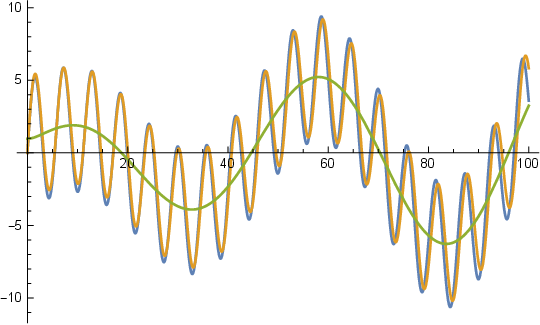}
	\caption{Plots of  the real (left panel)  and imaginary (right panel) part for the exact function $f(iy-2)$  (blue) and for the  approximations according to Eqs.~(\ref{f16}) (green) and (\ref{ha114}) (orange) for $a=3$.
	}
	\label{23}       
\end{figure}

\begin{figure}
	\includegraphics[width=0.49\textwidth]{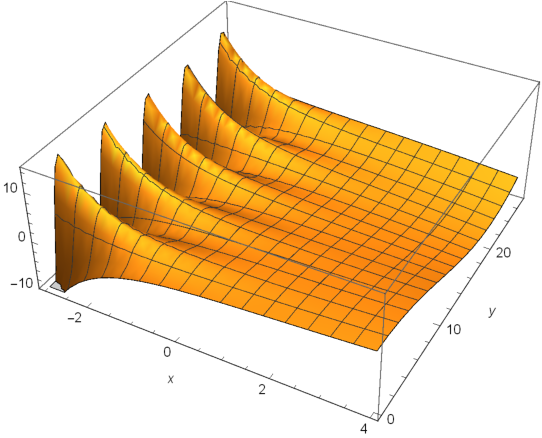}
		\includegraphics[width=0.49\textwidth]{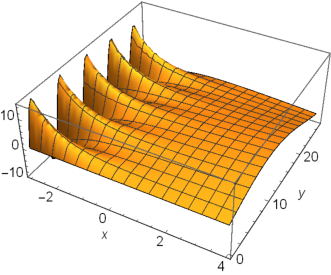}
	\caption{3D plots of  the real (left panel)  and imaginary (right) part of the function $f(x+iy,a)$ for $a=3$.
	}
	\label{12}       
\end{figure}

For $a>1$ the behavior of the numerical results becomes still more complicated and interesting: for Re$z<0$ the low-frequency oscillations become modulated by high-frequency ones (Figs.3, 4).
To clarify this behavior, we perform an analytical consideration.

The Taylor series expansion for the function (\ref{ha11}) is obtained from that for the error function
\begin{equation}
	{\rm erf} (z)=\frac {2}{\sqrt \pi}\sum_{n=0}^{\infty }(-1)^{n}\frac{z^{2n+1}}{n!(2n+1)}
\end{equation}%
and has the form
\begin{equation}
	\phi(z)=1.04\sum_{n=1}^{\infty }(-1)^{n-1}\frac{z^{n}}{(n-1)!(2n-1)},
	\label{h111}
\end{equation}%
being convergent for any complex $z$. Therefore, $\phi(t)$  has good analytical properties and can be safely continued into the complex plane.
Thus the numerical calculations for complex $z$ can be performed by  using the integration in (\ref{sa}), approximation (\ref{77}) and expansion (\ref{h111}). The results are in agreement, although a  small oscillation frequency shift  occurs. To make correction, we use an interpolation formula which takes into account this shift,
\begin{equation}
	f(t,a) =\left( \frac{2\pi t}{e a}\right) ^{1/2}
	{\rm erf}\left( \frac{ta}{e}\right) ^{1/2}	\exp \left(\frac{t}{e a}(1+ 1.047(a-4))^{-1} \right)
	\label{ha114}
\end{equation}
(note that for $a=4$ the shift is absent). The fitting is presented in Fig.~\ref{23}.

\begin{figure}
	\includegraphics[width=0.59\textwidth]{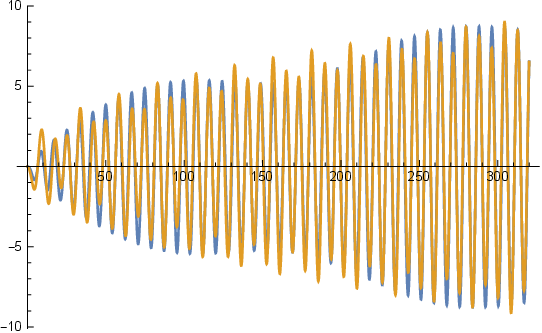}
	\caption{Plots of  the exact $f(iy - 2,a)$  (blue) and of the  approximation according to  (\ref{77}) (orange) for $a=0.48$.
	}
	\label{mod}       
\end{figure}

Thus $\phi(z)$  plays the role of a modulation factor.
To explain the origin of high-frequency oscillations, we can use the asymptotics at large $|z|$ \cite{error}
\begin{equation}
	{\rm erf}(z) \simeq 1 - \exp(-z^2)/ (\sqrt \pi z),
	\label{h11}
\end{equation}%
to obtain
\begin{equation}
	\phi(x+iy,a) \simeq \left( \frac{2\pi t}{e a}\right) ^{1/2}\left(1 - \frac {\exp (-a(x+iy)/e)}{\sqrt{\pi a(x+iy)/e}}\right)
	\label{h121}.
\end{equation}%
We see that the ratio of the high and low frequencies makes up $a^2$, and the amplitude of  high-frequency oscillations increases with increasing $-x>0$ (Fig.~6). For $x>0$, these oscillations are suppressed.

For $a=1+\delta$,  slow modulation occurs with a frequency $[(1+\delta) - 1/(1+\delta)]/e\simeq 2\delta/e$ due to an interference.

With decreasing $a$, this picture fails at $a \simeq 1/2$.
Again, in the crossover region below  $a=1/2$ the modulations become very slow, unlike the approximation (\ref{77}) (see Fig.~\ref{mod}).

For $a< 1/e$ no modulations take place, and only one frequency survives (remember that  the saddle point $y_0$ in (\ref{y}) is absent under this condition). As follows from (\ref{f16b}), the oscillations for small $a$  have the period of $2 \pi i/ \ln (1/a)$,
which passes into the period of $2 \pi ie/a$ with increasing $a$ through a maximum of $2 \pi ie$ for $a=1$.

The approximation (\ref{f16b}) works well on the imaginary axis and for Re~$z>0$, but does not properly describe the behavior at negative real part, as demonstrated by Fig.~8.

For negative $a$ we have
\begin{equation}
	f(t,a)=t\int_{0}^{1}(\cos ( \pi tx) +i\sin ( \pi tx))(-ax)^{-tx}dx,
	\label{h1213}
\end{equation}%
so that we have a strongly oscillating function  on the real axis for both the real and imaginary parts.
For $a=-1$  the oscillations have  the period of 2 (Fig.~9). For negative $a \neq -1$, additional real exponential factors occur.

\begin{figure}
	\includegraphics[width=0.49\textwidth]{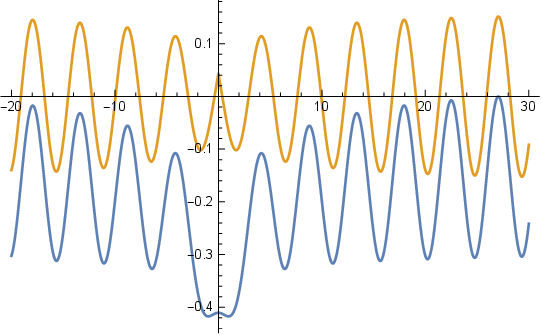}
		\includegraphics[width=0.49\textwidth]{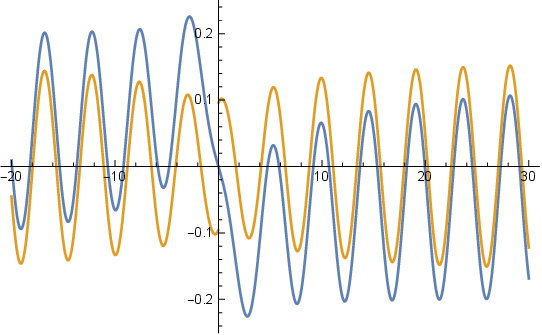}
	\caption{Plots of  the real (right panel) and imaginary (left panel)  parts of  the function $f(it-2,1/4)$, calculated numerically (blue) and according to the approximation (\ref{f16b}) (orange).
	}
	\label{-4}
\end{figure}

\begin{figure}
	\includegraphics[width=0.59\textwidth]{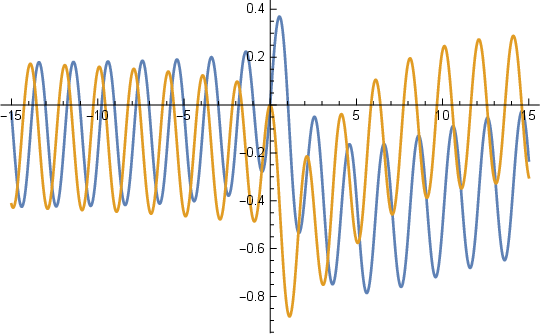}
	\caption{Plots of  the real (blue) and imaginary (orange) parts of  the function $f(t,-1)$.
	}
	\label{-1}       
\end{figure}

\begin{figure}
	\includegraphics[width=0.49\textwidth]{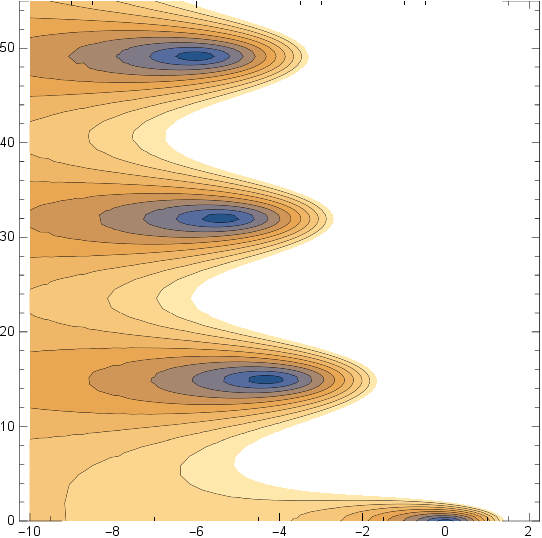}
		\includegraphics[width=0.49\textwidth]{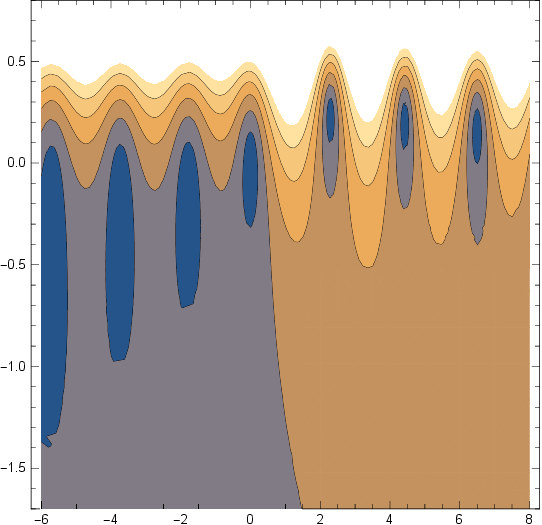}
	\caption{Contour plot of the function $|f(x+iy,a)|$ demonstrating occurrence of zeros for $a=1$ (left panel) and $a=-1$ (right panel).
	}
	\label{3}       
\end{figure}


\section{Zeros in the complex plane}

For positive $a$, the function $f(t,a)$
has non-trivial zeros lying
in the left complex half-plane  $x<0$.
Their position can be estimated from the zeros of the asymptotics for the error function ${\rm erf}(\sqrt t)$ which contains  according to (\ref{h11}) a periodic factor of $\exp (-t)$ with the period of $2\pi i$.

A more exact calculation performed in Ref.\cite{error}
gives the following values for the zeros of ${\rm erf} (s)$:
\begin{equation}
s_n = (1+i)\lambda_n  - \frac {1-i}{8 \lambda_n}\ln (2 \pi\lambda_n^2) +... \label{zero}
\end{equation}%
with $\lambda_n=\sqrt{(n-1/8)\pi}$.
Then, according to (\ref{haf11}), we can write down for the zeros of $f(x+iy)$
\begin{equation}
	t_n \simeq e s_n^2 \simeq 2\pi e (n-1/8)i - \frac {e}{2} \ln (2 \pi^2 (n-1/8)).
	  \label{zero1}
\end{equation}%
As demonstrate numerical calculations (Fig.~10, left panel), Eq.(\ref{zero1}) well describes the behavior of  the zeros of exact $f(x+iy)$.
It is interesting that we have a simple $n$-linear dependence of Im~$t_n$ in the complex plane, unlike the situation for the error function.

Similar results are obtained for $a > 1$ with Im~$t_n =  2\pi e (n-1/8)/a$, the poles being connected with high-frequency oscillations.
Moreover, $f(t,a)$ still has  zeros for small $a$, where the approximation (\ref{77}) fails. Indeed, since the periodic dependence of $f(t,a)$
along the imaginary axis  (with the  period of $2 \pi i /\ln a$)  retains, we have Im~$t_n \simeq  2\pi  (n-1/8)/ \ln a$.

In spirit of the Weierstrass theorem, we can write down in terms of the
zeros
\begin{equation}
f(t)=t\exp (h(t))\prod\limits_{n=1}^{\infty }\left\vert 1-\frac{t}{t_{n}}%
\right\vert ^{2}
\end{equation}%
where $h(t)$ is an entire function. In the simplest approximation $t_{n}\simeq 2\pi nei$ we have
\begin{equation}
f(t)\simeq t\exp (h(t))\prod\limits_{n=1}^{\infty }\left( 1+\left( \frac{t}{%
2\pi ne}\right) ^{2}\right)
\end{equation}%
Using Euler's infinite product representation%
\begin{equation*}
\sin x=x\prod\limits_{n=1}^{\infty }\left( 1-\frac{x^{2}}{\pi ^{2}n^{2}}%
\right)
\end{equation*}%
for $x=it/(2e)$ we obtain%
\begin{equation}
f(t)\simeq \exp (h(t))\sinh \frac{t}{2e}
\end{equation}%
Then from the asymptotics (\ref{f16}) we can put  $h(t)\simeq t/2e$  to estimate%
\begin{equation}
f(t)\ \propto \exp \frac{t}{e}-1
\end{equation}%
Although this approach is very crude, it may roughly describe the behavior
in the complex plane near the imaginary axis.

Qualitatively similar  pictures of oscillations and poles (but at different scales) occur for the function $\tilde{f}(t,a)$  (\ref{sad}).

For negative  $a$ the zeros reside slightly below the real axis. One can see in Fig.~10 (right panel) a correspondence with Fig.~9.



\section{Power-law and logarithmic generalizations}

 The function $f(t)$ can be generalized by including a power in the exponent, 
\begin{equation}
f^\mu(t)= t\int_{0}^{1}x^{-tx^\mu}dx.  \label{s334}
\end{equation}%
Similar to (\ref{s31})-(\ref{s311}), we can perform for $\mu>0$ the series expansion and integration to derive the relation
\begin{equation}
f^\mu(t)=\sum_{n=0}^{\infty }
\frac{t^{n+1}}{(\mu n+1)^{n+1}},\label{ffff}
\end{equation}
which is in agreement with (\ref{s16}), as is evident after changing the variable $z=x^\mu$ in (\ref{s334}).
In particular, for $t=1$ and $t=-1$ we obtain generalizations of Bernoulli's identities (\ref{s1}), (\ref{s2}).

For integer $m=1/\mu$ we obtain
\begin{equation}
f^{1/m}(t)=\sum_{n=m}^{\infty }
\frac{t^{n+1-m}m^n}{n^{n+1-m}} . \label{fff4}
\end{equation}

For large positive $t$, $f^\mu(t)$ exponentially increases with $t$,
\begin{equation}
	f^\mu(t)=\left( \frac{2\pi e t}{\mu}\right) ^{1/2}\exp \left(\frac{t}{e \mu} - \frac{1}{ \mu}\right).
	\label{f161}
\end{equation}
This asymptotics is derived similar to (\ref{f16}) with the saddle point being $y_0=1/\mu$.

For $t \rightarrow -\infty$, $f^\mu(t)$ tends to $-1$ passing through a minimum provided that $0.44\leq\mu<1$ and  without a minimum provided that $\mu\leq 0.43$, and rapidly decreases  provided that  $\mu>1$, a minimum being present at $1<\mu\leq 1.7$. 

At $\mu<0$,  $f^\mu(t)$ is determined only at negative $t$ and also tends to $-1$  for $t \rightarrow -\infty$.

The behavior of the function $f^\mu(t)$ in the complex plan is qualitatively similar to that of $f(t)$.

Another generalization of the function $f(t)$ can be written down as
\begin{equation}
f_l(t)=(-1)^l t\int_{0}^{1}x^{-tx}\ln^l(x)dx=(-1)^l t\int_{0}^{1}\ln^l(x^{x^{-tx/l}})dx.  \label{s333}
\end{equation}%
Similar to (\ref{s31})-(\ref{s311}), we can perform integration to derive the relation
\begin{equation}
f_l(t)=\sum_{n=1}^{\infty }\frac{(n+l-1)!}{n^{l}(n-1)!} 
\frac{t^{n}}{n^{n}}  \label{fff}
\end{equation}
which yields analogues of Bernoulli's identities (\ref{s1}), (\ref{s2}):
\begin{eqnarray}
\sum_{n=1}^{\infty }\frac{(n+l-1)!}{n^{l}(n-1)!} \frac{1}{n^{n}} &=&(-1)^l\int_{0}^{1}x^{-x}\ln^l(x)dx  \label{s111}, \\
\sum_{n=1}^{\infty }\frac{(n+l-1)!}{n^{l}(n-1)!} \frac{(-1)^{n-1}}{n^{n}} &=&(-1)^l\int_{0}^{1}x^{x}\ln^l(x)dx  .\label{s222}
\end{eqnarray}

Note that $f_1(t)=f_0(t)=f(t)$ identically, which can be also demonstrated by integration by parts in (\ref{s333}). Moreover, for $l>1$ we have $f_l(t)\simeq f(t)$ for large positive $t$, since large values of $n$ yield main contributions to the sum in (\ref{fff}). For $t \rightarrow -\infty$, $f_l(t)$ tends to a negative constant which rapidly increases with $l$, no minimum in $t$-dependence being present for $l>1$.


\begin{figure}
	\includegraphics[width=0.49\textwidth]{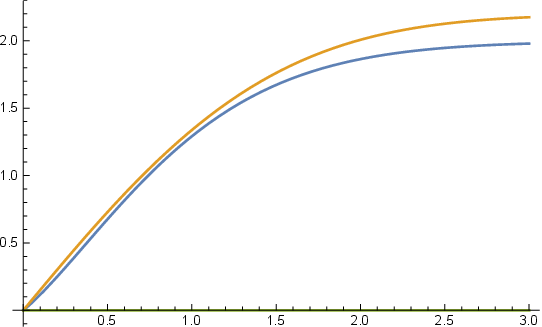}
	\includegraphics[width=0.49\textwidth]{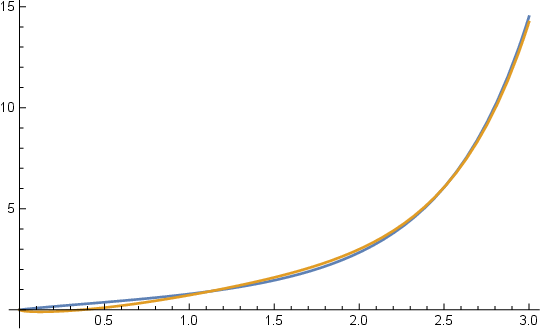}
	\caption{Left panel: Plots of  the exact $F(1,\lambda)$  (blue) and of the  approximation  Eq. (\ref{77}) (orange) according to Eq. (\ref{s15}).
		Right panel:	Plots of  the exact  $F(-1,\lambda)$ (blue) and of the  approximation (\ref{faaa}) (orange), the lines at large $\lambda$ being almost indistinguishable.
	}
	\label{9}      
\end{figure}

\section{Conclusions}
To conclude,  investigation of the function  $f(t,a)$ provides a rich picture including crossovers with changing the parameter $a$, which are connected with the presence of saddle points.
The consideration of Sect.~2 demonstrates usefulness of the half-derivative  concept in our case, which is a rather rare situation.
Using the error function provides a reliable approximation  describing non-trivial numerical results in the complex plane including modulated oscillations and poles in the left-half plane. Also asymptotics for negative $t$ are obtained with a high accuracy.

The integral with varying integration limits
is reduced to the function $f(t,a)$ by a linear variable change,
\begin{equation}
F(t,\lambda)=	t\int_{0}^{\lambda}x^{-tx}dx 	= \lambda f(\lambda t,\lambda) \label{s15}.
\end{equation}%
Examples, which demonstrate comparison of numerical integration with the approximations elaborated, are presented in Fig.~11.
Therefore, the results obtained above shed light on analytical integration of the function $x^{tx}$.
Further investigation of analytical properties of the function $f(t,a)$ and its generalizations would be of interest.


\begin{thebibliography}{99}

\bibitem{552}
Johann Bernoulli, Demonstratio methodi analyticæ, qua determinataest aliqua quadratura exponentialis per Seriem, Opera omnia, vol. 3 (1697).

\bibitem{553}
W. Dunham, The Calculus Gallery: Masterpieces from Newton to Lebesgue. Princeton, NJ: Princeton University Press,  2005.

\bibitem{Pi}
Gabor Roman, Extension of the Sophomore’s Dream, An. st.  Univ. "Ovidius" Constanta. Seria Matematica,
\textbf{29},  211 (2021).

\bibitem{Pi1}
Jean Jacquelin, Sophomore's Dream Function,  2010.
https://tetrationforum.org/attachment.php?aid=788

\bibitem{frac}
S. Samko, A. Kilbas, O. Marichev, Fractional Integrals and Derivatives: Theory and Applications,
CRC Press, 1993.

\bibitem{error}
Henry E. Fettis, James C. Caslin and Kenneth R. Cramer,
Complex Zeros of the Error Function and of
the Complementary Error Function,
Mathematics of Computation \textbf{27}, 401 (1973).

\end{thebibliography}
\end{document}